\magnification=1200

\tolerance=500

\font \twelvebf=cmbx12

\def\dem{\noindent {\bf D\'emonstration.}\enspace \nobreak }
\outer \def \th #1. #2\par{ \medbreak
\noindent {\bf#1. \enspace} {\sl#2 }\par
\ifdim \lastskip< \medskipamount \removelastskip \penalty55 \medskip \fi}
\def \expl#1 {\medbreak \noindent {\bf Exemple #1.}\enspace }

\def \mapdown#1{\Big\downarrow
\rlap{$\vcenter{\hbox{$\scriptstyle#1$}}$}} 
\def \rema#1 {\medbreak \noindent {\bf Remarque #1.}\enspace }
\def \remas#1 {\medbreak \noindent {\bf Remarques #1.}\enspace }

\def\tas {{\hbox{\bf .}}}

\def \titre#1{\medbreak \noindent {\bf #1.}\medbreak}

\def \cExt {{{\cal E}\it xt\,}}

\def \bZ {{\bf Z}}
\def \cB {{\cal B}}
\def \cE {{\cal E}}
\def \cL {{\cal L}}
\def \cP {{\cal P}}
\def \cF {{\cal F}}

\def \cJ {{\cal J}}
\def \cO {{\cal O}}
\def \cM {{\cal M}}
\def \cN {{\cal N}}
\def \cS {{\cal S}}

\centerline {\twelvebf  Minimalit\'e des courbes sous-canoniques}

\vskip 2 cm

\titre {0. Introduction}

Soient ${\bf P}^3$ l'espace projectif de dimension 3 sur un
 corps $k$ alg\'ebriquement clos et $R=k[X,Y,Z,T]$ l'anneau de polyn\^omes
associ\'e. Il y a des liens forts entre les  faisceaux coh\'erents 
(ou les fibr\'es) sur ${\bf P}^3$, les  $R$-modules gradu\'es de longueur finie
et les courbes localement Cohen-Macaulay de ${\bf P}^3$. Rappelons les
propri\'et\'es suivantes :

\vskip 0.3cm

L'\'equivalence stable est d\'efinie sur l'ensemble des fibr\'es :

\th {D\'efinition 0.1}.  Deux fibr\'es $\cF$ et $\cF'$ sur ${\bf P}^3$ sont dits
stablement isomorphes s'il existe des fibr\'es dissoci\'es
(c'est-\`a-dire sommes directes de faisceaux inversibles) $\cL$ et $\cL'$ et
un isomorphisme $\cF\oplus \cL\simeq \cF'\oplus \cL'$.

Dans [HMDP], nous avons d\'efini sur l'ensemble des faisceaux coh\'erents sur
${\bf P}^3$ la relation d'\'equivalence de
pseudo-isomorphisme  :

\th {D\'efinition 0.2}. Soient $\cN$ et $\cN'$  des faisceaux coh\'erents
sur ${\bf P}^3$  et soit
$f$ un morphisme de
$\cN$ dans $\cN'$. On dit que $f$ est un pseudo-isomorphisme   (en abr\'eg\'e
un {\it psi}) s'il induit~:{\hfil \break } 0) un isomorphisme  $H^0\cN(n) \to
H^0\cN'(n)$ pour tout $n \ll 0$,  {\hfil \break }
1) un isomorphisme  $H^1_*\cN \to H^1_*\cN'$ et {\hfil \break } 
2) une injection $H^2_*\cN \to H^2_*\cN'$.{\hfil \break } 
Deux faisceaux
coh\'erents seront dits pseudo-isomorphes s'il existe
une cha\^\i ne de {\it psi}\  qui les joint.

C'est une extension de  l'\'equivalence stable au sens suivant :

\th {Proposition 0.3}. L'application canonique de l'ensemble ${\cS
tab}$ des classes d'iso\-mor\-phisme stable de fibr\'es $\cF$ de
${\bf P}^3$ v\'erifiant $H^2_*\cF=0$ dans  l'ensemble ${\cP si}$ des
classes de pseudo-isomorphisme de faisceaux coh\'erents de dimension
projective $\leq 1$ est une bijection.

\dem L'injectivit\'e est une cons\'equence de [HMDP]2.11 et 2.8, la
surjectivit\'e de [HMDP]2.10.

Le lien entre les fibr\'es et les $R$-modules gradu\'es est le suivant :

\th {Proposition 0.4. (Horrocks, cf. [Ho])}. Soit $\cM_{f}$ l'ensemble des
classes d'iso\-mor\-phisme de
$R$-modules gradu\'es de longueur finie. L'application qui \`a un tel module
associe le faisceau associ\'e \`a son deuxi\`eme module de syzygies induit
une bijection de $\cM_{f}$ dans ${\cS tab}$, la bijection r\'eciproque
\'etant induite par l'application qui envoie un fibr\'e $\cF$ sur le module 
$H^1_*\cF$.

On en d\'eduit que l'application qui envoie un faisceau
$\cN$ de dimension projective $\leq 1$ sur le module  $H^1_*\cN$
induit une bijection de ${\cP si}$ sur $\cM_{f}$.

Passons maintenant aux courbes :

\th {Proposition 0.5. (Rao, cf. [R])}. Soit ${\cB il}$ l'ensemble des classes
de biliaison de courbes (localement Cohen-Macaulay) de ${\bf P}^3$. L'application
qui
\`a une courbe $C$ associe son module de Rao  $H^1_*\cJ_C$ induit une
bijection de ${\cB il}$ sur le quotient de $\cM_{f}$ par l'action de
d\'ecalage des degr\'es. 

\th {Corollaire 0.6}. L'application qui
\`a une courbe $C$ associe son faisceau d'id\'eaux $\cJ_C$ induit une
bijection de ${\cB il}$ sur le quotient de ${\cP si}$ par l'action de
tensorisation par une puissance du faisceau ${\cal O}_{\bf P}(1)$.

 Il y a deux mani\`eres
particuli\`eres de construire la bijection r\'eciproque ; en effet dans chaque
classe de ${\cP si}$, il y a des fibr\'es d'apr\`es 0.3, et des faisceaux
r\'eflexifs de rang 2 (cf. [MDP2]) :

-- soit $\cF$ un fibr\'e ; il
existe un faisceau dissoci\'e 
$\cP$, un entier $h$, une courbe $C$ et une suite exacte $0\to \cP \to \cF \to
\cJ_C(h)\to 0$. On associe \`a la classe de $\cF$ la classe
de biliaison de $C$ ;

-- soit  $\cN$ un faisceau r\'eflexif de rang 2 ; on lui  associe la
classe  de biliaison d'une courbe obtenue comme  sch\'ema des
z\'eros d'une section non nulle de $\cN(n)$ pour un entier $n$ bien choisi (en
particulier si $H^0\cN(n-1)=0$
et $H^0\cN(n)\neq 0$, toute section non nulle de $\cN(n)$ convient).

\vskip 0.3cm

Dans chaque classe de biliaison, il y a des  {\bf courbes minimales}, qui
r\'ealisent le plus petit d\'ecalage [cf. Mi], et qui permettent de d\'ecrire
toutes les courbes de la classe (cf. [MDP1]~V, [BBM]).

Dans chaque classe de ${\cP si}$, il y a,
 parmi les faisceaux r\'eflexifs de rang 2, des {\bf faisceaux
r\'eflexifs minimaux}, dont la troisi\`eme classe de Chern est minimale (cf.
[B]). 

Il est naturel de se demander s'il y a une relation entre les courbes
minimales et les faisceaux
r\'eflexifs minimaux (\`a d\'ecalage pr\`es), ce qui conduit \`a la  
question suivante~:

 \th {Question I}. Dans une classe de biliaison, la courbe minimale
est-elle section d'un faisceau r\'eflexif ?  

La r\'eponse 
 est n\'egative, comme on le voit facilement sur un
contre-exemple (cf. [B]). Si $M$ est un module de Koszul 
de type $(n_1,n_2,n_3,n_4)$,
c'est-\`a-dire un  quotient de $R$ par une suite r\'eguli\`ere
$(f_1,f_2,f_3,f_4)$ o\`u $f_i$ est de degr\'e $n_i$ avec
$n_1\leq n_2 < n_3\leq n_4$, la courbe minimale
associ\'ee au module $M\oplus M$ n'est pas une section d'un
faisceau r\'eflexif minimal.

\vskip 0.3cm

Dans une classe de ${\cP si}$, il n'y a pas toujours de fibr\'e de rang 2.
Lorsqu'il y en a, ce sont les \'el\'ements minimaux. Les sch\'ema des
z\'eros des sections de ces fibr\'es, lorsqu'ils sont de dimension 1, sont
des courbes {\bf sous-canoniques}. On peut alors poser la question suivante :

 \th {Question II}. Si une classe de biliaison contient des
courbes sous-canoniques,  la courbe minimale est-elle aussi 
sous-canonique ?

Par exemple, dans la classe de biliaison associ\'ee \`a un module de
Kozsul de type $(n_1,n_2,n_3,n_4)$
avec
$n_1\leq n_2 \leq n_3\leq n_4$ et $n_1+n_4=n_2+n_3$, la courbe minimale est
sous-canonique.

\vskip 0.3cm

Comme le montre A. Buraggina (cf. [B] 5), cette question est \'equivalente
\`a la question suivante, qui nous a \'et\'e pos\'ee par  Hartshorne et Ellia :

 \th {Question III}. Soit
$\cE$ un fibr\'e
 de rang 2 sur ${\bf P}^3$, $n$ un entier relatif tel que $H^0\cE(n-1)=0$
et $H^0\cE(n)\neq 0$, soit $C$ une courbe sch\'ema des z\'eros
d'une section non nulle de $\cE(n)$. Est-elle minimale dans sa classe
de biliaison ?

Dans cet article, nous donnons une r\'eponse positive aux questions II et
III (th\'eor\`eme 2.5).

 Dans le premier paragraphe, nous \'etudions, pour toute courbe $C$ trac\'ee
 sur une surface $Q$, le faisceau 
${\cal H}om_{{\cal O}_Q}({\cal J}_{C/Q},{\cal O}_Q)$ dont les sections
globales  sont li\'ees aux biliaisons \'el\'ementaires que l'on peut faire
\`a partir de la courbe (cf. 1.2) et \`a ses propri\'et\'es de mini\-malit\'e.
En particulier nous caract\'erisons les homomorphismes  non nuls et non
injectifs de ${\cal J}_{C/Q}$ dans ${\cal O}_Q(h)$ (cf. 1.8).

Le deuxi\`eme paragraphe est consacr\'e \`a la preuve du r\'esultat. La
m\'ethode est la suivante :  si $C$ est une courbe sous-canonique minimale
pour un fibr\'e pour laquelle on peut faire une biliaison
\'el\'ementaire descendante, il existe un entier $n<0$ et une section non
nulle de
${\cal O}_C(n)$. L'\'etude de la courbe contenue dans $C$ sur laquelle cette
section s'annule conduit \`a une contradiction.

{\medbreak \noindent {\bf Notations.}}
On d\'esigne par $k$ un corps alg\'ebriquement clos et par $R$
l'anneau de polyn\^omes 
$k[X,Y,Z,T]$.  L'espace projectif ${\bf P}^3_k$ sera not\'e simplement 
$ {\bf P}^3$ et son
faisceau structural  ${\cal O}_{\bf P}$. Si ${\cal F}$ est un ${\cal O}_{\bf
P}$-module nous  noterons $h^i \cF$ la dimension de l'espace
vectoriel $H^i \cF$, et $H^i_* \cF$ le $R$-module gradu\'e $ \bigoplus _{n\in
\bZ}H^i
\cF(n)$.

\medbreak
Une  courbe $C$ de $ {\bf P}^3$  est un sous-sch\'ema ferm\'e purement
de dimension 1, localement Cohen-Macaulay, d\'efini par un faisceau
d'id\'eaux $\cJ_C$. Son
faisceau dualisant est le faisceau  $\omega_C=\cExt^2_{{\cal O}_P}({\cal
O}_C,\omega_P)=\cExt^2_P({\cal O}_C,{\cal O}_{\bf P})(-4)$.

On note $$e(C)=  \sup \{ n \in {{\bf Z}}\;\mid\; h^1\cO_C(n) \neq 0 \; \},$$
$$s_0(C)=  \inf \{ n \in {{\bf Z}}\;\mid\; h^0\cJ_C(n) \neq 0 \; \}.$$

 Le module de Rao de $C$ : $M_C =H^1_*\cJ_C$ est un $R$-module gradu\'e de
longueur finie qui joue un  r\^ole important dans la classification des
courbes gauches.

 \titre {1. Etude du dual de l'id\'eal d'une courbe trac\'ee sur une surface}

Dans tout ce paragraphe, on d\'esignera par $Q$ une surface  de degr\'e
$s$, non n\'ecessairement int\`egre, de $ {\bf P}^3$ et par $q$ son
\'equation. Pour toute courbe $C$ trac\'ee sur $Q$, d\'efinie par un faisceau
d'id\'eaux $\cJ_C$, on va \'etudier le faisceau de ${\cal O}_Q$-modules
${\cal H}om_{{\cal O}_Q}({\cal J}_{C/Q},{\cal O}_Q)$. Les sections
globales de ce faisceau sont li\'ees aux propri\'et\'es de minimalit\'e de la
courbe, comme nous le rappelons ci-dessous.

\th {D\'efinition 1.1}. Une courbe $C$ est minimale dans sa classe
de biliaison si son module de Rao a le d\'ecalage minimum,
c'est-\`a-dire si pour toute courbe
$C'$ de la classe de biliaison de $C$ on a $M_{C'}\simeq M_C(-h)$ avec
$h\geq 0$.

Dans la description des classes de biliaison, on utilise 
l'op\'eration de biliaison \'el\'emen\-tai\-re, qui   s'obtient en 
pratiquant deux liaisons successives, l'une des surfaces liantes \'etant
commune aux deux liaisons. Plus pr\'ecis\'ement, on a le r\'esultat suivant :
  
\th { Proposition 1.2}. Soient $C$ et $C'$ deux courbes trac\'ees sur
$Q$ et soit $h \in \bZ$. Les conditions suivantes sont \'equivalentes~: {\hfil \break }
1) $C'$ est obtenue \`a partir de $C$ par une double liaison par des 
surfaces $(Q,S)$ et
$(Q,S')$ avec $\deg S' - \deg S = h$. {\hfil \break }
2) Il existe un homomorphisme injectif $u : {\cal J}_{C/Q}(-h) \to
{\cal O}_Q$  d'image ${\cal J}_{C'/Q}$. {\hfil \break }
On a alors $M_C\simeq M_{C'}(h)$.{\hfil \break }
On dit que $C'$ est obtenue \`a partir de $C$ par une {\bf biliaison
\'el\'ementaire
 de hauteur
$h$ sur $Q$}, ascendante  (resp.~ descendante) si $h>0$ (resp.~$h<0$).

\dem Voir [MDP1] III.2.3.

\th { Proposition 1.3. Le diagramme fondamental}. Soient $C$ une courbe
trac\'ee sur $Q$, $h\in\bZ$ et $u :  {\cal J}_{C/Q}(-h) \to{\cal O}_Q$ un
homomorphisme. On a un diagramme commutatif de suites exactes de ${\cal
O}_Q$-modules :{\hfil \break }
$$\matrix {0&\to &{\cal J}_{C/Q}(-h) &\to{j} &{\cal O}_Q(-h) &
\to{}&{\cal O}_C(-h)&\to & 0\cr
&&\mapdown{u}&&\mapdown{u^{\vee}}&&\mapdown{\theta}&& \cr 0&\to & {\cal
O}_Q &
\to{j^{\vee}}&{\cal H}om_{{\cal O}_Q}({\cal J}_{C/Q},{\cal O}_Q) & \to{}&
\omega_C(4-s) &\to & 0\cr}$$
o\`u $j : {\cal J}_{C/Q}\to{\cal O}_Q$ est l'injection canonique.

\dem Partant de la suite exacte :

$$0\to {\cal J}_{C/Q} \to{j} {\cal O}_Q 
\to{}{\cal O}_C\to  0$$
on obtient la premi\`ere ligne en la tensorisant par ${\cal O}_Q(-h)$ et la
deuxi\`eme ligne en lui appliquant le foncteur ${\cal
H}om_{{\cal O}_Q}(\tas,{\cal O}_Q)$. En effet, on a ${\cal
H}om_{{\cal O}_Q}({\cal O}_C,{\cal O}_Q)=0$ et
$\omega_C={\cal E}xt^1_{{\cal O}_Q}({\cal
O}_C,\omega_Q)={\cal E}xt^1_{{\cal O}_Q}({\cal O}_C,{\cal O}_Q)(s-4)$.

L'\'egalit\'e $u^{\vee}j=j^{\vee}u$, qui entra\^\i ne l'existence de
$\theta$, est une cons\'equence du lemme facile d'alg\`ebre suivant :

\th {Lemme 1.4}. Soit $A$ un anneau commutatif, $J$ un id\'eal de $A$, $j$
l'injection canonique de $J$ dans $A$ et
$u$ un homomorphisme $A$-lin\'eaire de $J$ dans $A$. Alors on a
$u^{\vee}j=j^{\vee}u$.

\dem  Soient $\alpha$ et $\beta$ deux \'el\'ements de $J$. On a
:

$u^{\vee}j(\alpha)(\beta)=j(\alpha)u(\beta)=\alpha
u(\beta)=u(\alpha \beta)$,
$j^{\vee}u(\alpha)(\beta)=u(\alpha)j(\beta)=\beta u(\alpha)=u(\alpha
\beta)$.

\vskip 0.3cm
Gardant les notations de 1.3, on en d\'eduit les deux r\'esultats suivants,
qui seront  utiles dans la suite :

\th {Corollaire 1.5}. Pour tout entier n\'egatif $h$, on a un isomorphisme 
$Hom({\cal J}_{C/Q},{\cal O}_Q(h))\simeq H^0\omega_C(4-s+h)$ qui \`a $u$
associe $\theta$.

\dem Cela r\'esulte de la suite exacte :
$$0\to H^0 {\cal O}_Q(h)\to Hom({\cal J}_{C/Q},{\cal O}_Q(h))\to
H^0\omega_C(4-s+h)\to 0$$

\th {Corollaire 1.6}. Soit $C'$ une courbe contenue dans $C$. Alors $u$ se
prolonge \`a ${\cal J}_{C'/Q}(-h) $ si et seulement si $\theta$ se factorise
\`a travers la projection ${\cal O}_C(-h)\to{\cal O}_{C'}(-h)$.

\dem Soient $j' : {\cal J}_{C'/Q} \to{\cal O}_Q$ et $i : {\cal
J}_{C/Q} \to {\cal J}_{C'/Q}$ les injections canoniques. Supposons
que $u$ se
prolonge en $u': {\cal J}_{C'/Q}(-h)\to {\cal O}_Q $. On a un diagramme
commutatif de suites exactes :
$$\matrix {0&\to &{\cal J}_{C/Q}(-h) &\to{j} &{\cal O}_Q(-h) &
\to{}&{\cal O}_C(-h)&\to & 0\cr
&&\mapdown{i}&&\parallel&&\mapdown{}&& \cr
0&\to &{\cal J}_{C'/Q}(-h) &\to{j'} &{\cal O}_Q(-h) &
\to{}&{\cal O}_{C'}(-h)&\to & 0\cr
&&\mapdown{u'}&&\mapdown{u^{'\vee}}&&\mapdown{\theta'}&& \cr
 0&\to & {\cal
O}_Q &
\to{j^{'\vee}}&{\cal H}om_{{\cal O}_Q}({\cal J}_{C'/Q},{\cal O}_Q) & \to{}&
\omega_{C'}(4-s) &\to & 0\cr
&&\parallel&&\mapdown{i^{\vee}}&&\mapdown{}&& \cr
 0&\to & {\cal
O}_Q &
\to{j^{\vee}}&{\cal H}om_{{\cal O}_Q}({\cal J}_{C/Q},{\cal O}_Q) & \to{}&
\omega_{C}(4-s) &\to & 0\cr}$$
L'\'egalit\'e $i^{\vee}u^{'\vee}=u^{\vee}$ montre que la compos\'ee des trois
fl\`eches verticales de droite n'est autre que $\theta$, qui se factorise
comme annonc\'e.

Inversement, si $\theta$ se factorise \`a travers la projection ${\cal
O}_C(-h)\to{\cal O}_{C'}(-h)$,  le diagramme pr\'ec\'edent dans lequel
on supprime la troisi\`eme ligne nous donne l'existence de la fl\`eche $u'$.

\vskip 0.3cm
D'apr\`es 1.2, l'\'etude de la minimalit\'e d'une courbe trac\'ee sur $Q$ est li\'ee \`a
l'existence d'homomorphismes {\bf injectifs} (et non surjectifs) de ${\cal
J}_{C/Q}$ dans
${\cal O}_Q(h)$ avec
$h$ n\'egatif. D'apr\`es 1.5, l'existence d'homomorphismes {\bf non nuls} 
de ${\cal J}_{C/Q}$ dans ${\cal O}_Q(h)$ est \'equivalente \`a l'existence de
sections non nulles du faisceau $\omega_C(4-s+h)$, donc, puisqu'on a \break
$H^0\omega_C(4-s+h)=H^1{\cal O}_C(s-4-h)$,  \`a l'in\'egalit\'e $h\geq
s-4-e(C)$ ; il sera possible d'en obtenir avec $h<0$ si et seulement si
$s_0(C)\leq s<e(C)+4$. Il  faut ensuite \'etudier quels sont les
homomorphismes non nuls et non injectifs, qui ne peuvent exister que si $Q$
n'est pas int\`egre.

\rema {1.7} Si $Q$ n'est pas int\`egre, posons $q=q_1q_2$ et soient  $Q_1$ et
$Q_2$ les surfaces correspondantes, $s_1$ et $s_2$ leurs degr\'es. Pour tout
$(i,j)\in \{ (1,2),(2,1)\} $ on a une suite exacte :
$$0\to {\cal O}_{Q_i}(-s_j) \to{\lambda_j}{\cal O}_{Q}\to{p_j}
{\cal O}_{Q_j}\to 0$$
o\`u $\lambda_jp_i$ est \'egal \`a la multiplication par $q_j:{\cal
O}_{Q}(-s_j)\to {\cal O}_{Q}$.

\th {Proposition 1.8}. Soient $C$ une courbe
trac\'ee sur $Q$, $h\in\bZ$,  $u :  {\cal J}_{C/Q}(-h) \to{\cal O}_Q$ un
homomorphisme non nul et $\theta$ la section de $\omega_C(4-s+h)$ qui lui
correspond. Les conditions suivantes sont \'equivalentes~:
{\hfil \break } 
i) $u$ n'est pas injectif, {\hfil \break }
ii) il existe une d\'ecomposition $q=q_1q_2$, o\`u $q_1$ et $q_2$ ne sont pas
constants, telle que
$q_1u=0$,{\hfil \break } 
iii) il existe une d\'ecomposition $q=q_1q_2$, o\`u $q_1$ et $q_2$ ne sont pas
constants, telle que, avec les notations de la remarque 1.6,  $u$ se factorise
par $\lambda_2 :{\cal O}_{Q_1}(-s_2) \to{\cal O}_{Q}$.{\hfil \break }
De plus, si $h<0$, elles sont encore \'equivalentes \`a :{\hfil \break }
iv) il existe une d\'ecomposition $q=q_1q_2$, o\`u $q_1$ et $q_2$ ne sont pas
constants, telle que
$q_1\theta=0$.

\dem  $i \Rightarrow ii$ : 
si $u$ n'est pas injectif, il en est de m\^eme de l'homomorphisme de
modules $ :I_C/(q)(-h) \to R/(q)$  associ\'e, qu'on d\'esignera encore par
$u$. Soient
$g$ un
\'el\'ement de
$I_C$ dont l'image $\overline g$ est un \'el\'ement non nul du
noyau de $u$ et
$q_1$ le pgcd de $q$ et $g$, de sorte qu'on a $q=q_1q_2$ et $g=q_1g'$, o\`u
$q_2$ et
$g'$ sont premiers entre eux. Puisque $\overline g$ n'est pas nul, $q_1$
est un diviseur strict de $q$, et $q_2$ n'est pas une constante. Pour
tout
$f$ dans
$I_C$ on a 
$gu(\overline f)=fu(\overline g)=0$ ; on en d\'eduit que si $f'$ rel\`eve
$u(\overline f)$,
$q=q_1q_2$ divise $gf'=q_1g'f'$, donc $q_2$ divise $g'f'$, $q_2$ divise $f'$,
et
$q$ divise
$q_1f'$. On a donc montr\'e que $q_1u=0$. Puisque $u$ n'est pas nul, ceci
prouve aussi que
$q_1$ n'est pas une constante.

 $ii \Leftrightarrow iii$ : puisque $\lambda_1$ est injectif,
$q_1u=\lambda_1p_2u$ est nul si et seulement si $p_2u$ est nul, ce qui
\'equivaut au fait que $u$ se factorise par $\lambda_2$.

 $iii \Rightarrow i$ : un homomorphisme ${\cal J}_{C/Q}(-h) \to{\cal
O}_{Q_1}(-s_2)$ ne peut pas \^etre injectif, car le support sch\'ematique de 
${\cal J}_{C/Q}$ contient strictement celui de ${\cal
O}_{Q_1}$.

$iv \Rightarrow ii$ si $h<0$ : si $q_1\theta=0$, $q_1u :  {\cal J}_{C/Q}(-h)
\to{\cal O}_Q\to {\cal O}_Q(s_1)$ se prolonge \`a ${\cal O}_Q(-h)$, autrement
dit il existe $v :{\cal O}_Q(-h)\to {\cal O}_Q(s_1)$ tel qu'on ait $q_1u=vj$.
On a alors $p_2vj=0$, donc $p_2v$ se factorise par la projection ${\cal
O}_Q(-h)\to {\cal O}_{C}(-h)$ compos\'ee avec un homomorphisme
${\cal O}_{C}(-h)\to {\cal O}_Q$ qui est nul pour des raisons
de profondeur. Alors $v$ se factorise \'egalement par $\lambda_1 :{\cal
O}_{Q_2} \to{\cal O}_{Q}(s_1)$, donc il existe  $w :{\cal O}_Q(-h)\to {\cal
O}_{Q_2}$ tel qu'on ait $v=\lambda_1w$. Puisque $h$ est n\'egatif, $w$ est
nul d'o\`u le r\'esultat.

\remas {1.9} 

1) Si on ne suppose plus que $u$ n'est pas nul, les conditions ii)
 iii) et iv) restent valables, \`a condition de supposer seulement que $q_2$
n'est pas  constant, c'est-\`a -dire que $q_1$ est un diviseur strict  de
$q$.

2) Si $h<0$ et si on a  $q=q_1q_2$, o\`u $q_1$ et $q_2$ ne sont pas
constants et si la  surface $Q_1$ d'\'equation $q_1$ contient $C$, $q_1$
annule
$H^0\omega_C(4-s+h)$ donc il n'existe pas d'homomorphisme injectif  $u : 
{\cal J}_{C/Q}(-h) \to{\cal O}_Q$.

\th {Corollaire 1.10}. Soient $C$ et $C'$ deux courbes trac\'ees sur
$Q$ telles que $C'$ soit contenue dans $C$, $h \in \bZ$. Soit  $u' :  {\cal
J}_{C'/Q}(-h) \to{\cal O}_Q$ et  $u :  {\cal J}_{C/Q}(-h) \to{\cal O}_Q$
sa restriction. Alors {\hfil \break }
1)  si $u'$ n'est pas nul, il en est de m\^eme de $u$ ;{\hfil \break }
2) si $u$ est injectif, il en est de m\^eme de $u'$.

\dem Si $u$ est nul, $u'$ se factorise par  par la projection ${\cal
J}_{C'}(-h)\to {\cal J}_{C'}/{\cal J}_{C}(-h)$ compos\'ee avec
un homomorphisme
${\cal J}_{C'}/{\cal J}_{C}(-h)\to {\cal O}_Q$ qui est nul pour des raisons
de profondeur, d'o\`u 1).

Si $u'$ n'est pas injectif, d'apr\`es 1.7, il existe une d\'ecomposition
$q=q_1q_2$, o\`u 
$q_1$ et $q_2$ ne sont pas constants, telle que $q_1u'=0$, mais alors on a
aussi par restriction $q_1u=0$, donc $u$ n'est pas injectif.

\th {Corollaire 1.11}. Soit $C$ une courbe. On ne peut pas faire \`a partir
de 
$C$ de biliaison \'el\'ementaire 
de hauteur n\'egative si et seulement si pour tout $s\geq s_0(C)$, pour tout
$h<0$, pour toute surface $Q$ de degr\'e $s$ d'\'equation $q$ contenant $C$,
il existe un diviseur strict $q_1$ de $q$ qui annule $H^0\omega_C(4-s+h)$.

\dem D'apr\`es 1.8, il suffit de voir que si toute section de
$\omega_C(4-s+h)$ est annul\'ee par un diviseur strict  de $q$, il en existe
un qui les annule toutes. Cela provient du fait que si un espace vectoriel
est  r\'eunion d'un nombre fini de sous-espaces vectoriels, il est \'egal \`a
l'un d'entre eux.

\expl {1.12} Consid\'erons un module de Koszul, c'est-\`a-dire un module
quotient de $R$ par une suite r\'eguli\`ere
$(f_1,f_2,f_3,f_4)$ o\`u $f_i$ est de degr\'e $n_i$
avec
$n_1\leq n_2 \leq n_3\leq n_4$. On pose $\mu = \sup \,(n_1+n_4 , n_2+n_3)$.
Toute courbe minimale associ\'ee a son id\'eal gradu\'e engendr\'e par
des polyn\^omes (rang\'es par degr\'es croissants) $ff_1^2$, $f_1f_2$,
$gf_2^2$,
$ff_1f_4+gf_2f_3$,  o\`u $f$ et $g$ sont 
des polyn\^omes
homog\`enes  de degr\'es respectifs $\mu -n_1-n_4$ et $
\mu-n_2-n_3$, non nuls et  tels que $f$, $g,$ et les $f_i$ soient deux \`a
deux sans facteur commun (cf. [MDP1] IV 6). 

On a $s_0(C) = \mu+n_1-n_4 $ et $ e(C) = 2\mu -n_3-n_4-4$. Les valeurs de $s$
et $h$ \`a consid\'erer sont celles qui v\'erifient $s_0(C)\leq s\leq 
e(C)+4+h$, ou encore  $\mu+n_1-n_4 \leq s\leq  2\mu -n_3-n_4+h$. Quand il en
existe (c'est-\`a-dire si on n'a pas \`a la fois $n_1=n_2$, $n_3=n_4$), les
\'equations des surfaces $Q$ correspondantes sont dans l'id\'eal 
$(ff_1^2,f_1f_2)$, donc sont toutes divisibles strictement par $f_1$, et on
v\'erifie que les sections de $\omega_C(4-s+h)$ sont annul\'ees par $f_1$.

\titre {2. Minimalit\'e des courbes sous-canoniques}

\th {D\'efinition 2.1}. Une courbe est dite sous-canonique s'il existe  un
fibr\'e
$\cE$ de rang 2 sur ${\bf P}^3$, un entier relatif $n$ et une
section non nulle de $\cE(n)$ dont le sch\'ema des z\'eros est $C$. On a
alors $\omega_C\simeq {\cal O}_C(2n+c_1-4)$, o\`u $c_1$ est la premi\`ere classe
de Chern de $\cE$. Une telle courbe est dite minimale pour
$\cE$  si
$\cE(n-1)$ n'a pas de section globale non  nulle.

\th {Proposition 2.2}. Soit $C$ une courbe sous-canonique, et
$C'$ une courbe contenue dans $C$ distincte de $C$. Alors il existe une
courbe $C''$ contenue dans $C$ et des isomorphismes : ${\cal J}_{C'}/{\cal
J}_C\simeq \omega_{C''}(-\alpha)$, ${\cal J}_{C''}/{\cal
J}_C\simeq \omega_{C'}(-\alpha)$ o\`u $\alpha$ est  l'entier qui v\'erifie
$\omega_C\simeq {\cal O}_C(\alpha)$. De plus, ${\cal J}_{C''}$ (resp.~ ${\cal
J}_{C'}$) est l'annulateur de ${\cal J}_{C'}/{\cal J}_C$ (resp.~ ${\cal
J}_{C''}/{\cal J}_C)$.

\dem  On consid\`ere la suite exacte $0\to {\cal J}_{C'}/{\cal J}_C\to {\cal
O}_C\to {\cal O}_{C'}\to 0$ et on lui applique le foncteur  ${\cal
H}om_{{\cal O}_{\bf P}}(\tas,{\cal O}_{\bf P})$. Le support de ${\cal
J}_{C'}/{\cal J}_C$ \'etant de dimension 1, le faisceau  ${\cal E}xt^1({\cal
J}_{C'}/{\cal J}_C,{\cal O}_{\bf P})$ est nul. On obtient la suite
exacte :
$$0\to {\cal E}xt^2({\cal O}_{C'},{\cal O}_{\bf P})\to {\cal E}xt^2({\cal
O}_{C},{\cal O}_{\bf P})\to {\cal E}xt^2({\cal J}_{C'}/{\cal J}_C,{\cal
O}_{\bf P})\to 0$$ 
Sachant qu'on a des isomorphismes :  ${\cal E}xt^2({\cal
O}_{C'},{\cal O}_{\bf P})\simeq \omega_{C'}(4)$, ${\cal E}xt^2({\cal
O}_{C},{\cal O}_{\bf P})\simeq \omega_{C}(4)\simeq {\cal O}_C(\alpha+4)$, on en
d\'eduit un isomorphisme : ${\cal E}xt^2({\cal J}_{C'}/{\cal J}_C,{\cal
O}_{\bf P})\simeq {\cal O}_{C''}(\alpha+4)$ o\`u $C''$ est un sous-sch\'ema
ferm\'e de
$C$ qui v\'erifie donc ${\cal J}_{C''}/{\cal
J}_C\simeq \omega_{C'}(-\alpha)$.

D'autre part, en appliquant de nouveau le foncteur ${\cal
H}om_{{\cal O}_{\bf P}}(\tas,{\cal O}_{\bf P})$ \`a la suite exacte obtenue :
$$0\to \omega_{C'}(4)\to \omega_{C}(4)\to {\cal O}_{C''}(\alpha+4)\to 0$$ 
et en
tenant compte des isomorphismes canoniques ${\cal E}xt^2(\omega_{C}(4),{\cal
O}_{\bf P})\simeq {\cal O}_{C}$ et \break ${\cal E}xt^2(\omega_{C'}(4),{\cal
O}_{\bf P})\simeq {\cal O}_{C'}$ on obtient  la suite exacte :
$$0\to {\cal E}xt^2({\cal O}_{C''}(\alpha+4),{\cal O}_{\bf P})\to {\cal O}_C\to
{\cal O}_{C'}\to {\cal E}xt^3({\cal O}_{C''}(\alpha+4),{\cal O}_{\bf P})\to 0$$
 ce qui prouve
que $C''$ n'est pas vide, que ${\cal E}xt^3({\cal O}_{C''},{\cal O}_{\bf
P})=0$, donc que
$C''$ est une courbe (localement Cohen-Macaulay) et que ${\cal E}xt^2({\cal
O}_{C''}(\alpha+4)\simeq \omega_{C''}(-\alpha)\simeq {\cal J}_{C'}/{\cal J}_C$.

La deuxi\`eme assertion r\'esulte du fait que l'annulateur de $\omega_{C'}$ 
(resp.~
$\omega_{C''})$ n'est autre que ${\cal J}_{C'}$ (resp.~ ${\cal J}_{C''}$).

\th {Corollaire 2.3}. Soient $C$ une courbe sous-canonique,  $Q$ une surface de degr\'e
$s$ contenant $C$, $h$ un entier, $u :  {\cal J}_{C/Q}(-h) \to{\cal O}_Q$ un
homomorphisme non nul et $\theta$ la section de $\omega_C(4+h-s)=
{\cal O}_C(\alpha+4+h-s)$ correspondante. Il
existe deux courbes $C'$ et $C''$ contenues dans $C$, une suite exacte $0\to
\omega_{C'}(-\alpha)\to {\cal O}_{C}\to {\cal O}_{C''}\to 0$ telles que  $C''$
soit la plus grande courbe (\'eventuellement vide) contenue
 dans le support du conoyau de $\theta$, que
${\cal J}_{C'}/{\cal J}_C$ soit le noyau de $\theta$, que $u$ se prolonge
\`a ${\cal J}_{C'/Q}(-h)$ et que ${\cal J}_{C''}$ soit l'annulateur de
${\cal J}_{C'}/{\cal J}_C$.{\hfil \break }
Si $h<(s-\alpha-4$, $C''$ n'est pas vide.

\dem L'image de $\theta :{\cal O}_{C}\to \omega_C(4+h-s)={\cal
O}_{C}(\alpha+4+h-s)$ est un quotient de ${\cal O}_{C}$, autrement dit
$\theta$ peut se factoriser de la mani\`ere suivante :
$\theta =\epsilon p$, o\`u
$p$ est la projection de ${\cal O}_{C}$ sur un quotient ${\cal
O}_{C'}$ et
$\epsilon $ est une injection ${\cal O}_{C'}\to {\cal O}_{C}(\alpha+4+h-s)$.
Ceci entra\^\i ne en particulier que tous les associ\'es de ${\cal O}_{C'}$
sont de dimension 1, donc que
$C'$ est une courbe (non vide car $\theta\neq 0$).

Le conoyau de $\theta$ est de la forme ${\cal O}_{Z}(\alpha+4+h-s)$ o\`u $Z$
est un sous-sch\'ema ferm\'e de $C$. Soit $C''$ la plus grande courbe contenue
 dans $Z$, qui
est \'egale \`a $Z$ en dehors d'un nombre fini de points. D'apr\`es 2.2, si
$C''$ n'est pas vide, il existe une courbe
$C'_1$ et une suite exacte
$0\to
\omega_{C'_1}(-\alpha)\to {\cal O}_{C}\to {\cal O}_{C''}\to 0$, et $\epsilon $
se factorise par une injection ${\cal O}_{C'}\to \omega_{C'_1}(4+h-s)$ qui
est un isomorphisme en-dehors d'un nombre fini de points. On en d\'eduit que
$C'$ et $C'_1$ sont \'egales.

Le fait que $u$ se prolonge \`a ${\cal J}_{C'/Q}(-h) $ r\'esulte de 1.6.

On v\'erifie que l'assertion est encore vraie (mais sans int\'er\^et) si
$C''$ est vide, ce qui correspond au cas o\`u $\theta$ est injective.

Si $h<(s-\alpha-4$, $\alpha+4+h-s$ est strictement n\'egatif donc $\theta$ n'est
pas injective (sinon son conoyau serait de longueur finie et aurait une
caract\'eristique de Hilbert strictement n\'egative).

 \th {Proposition  2.4}. Soient
$\cE$ un fibr\'e
 de rang 2 sur ${\bf P}^3$,  $C$ une courbe sous-canonique minimale pour
$\cE$,  $Q$ une surface  de degr\'e $s$ contenant $C$ et $h$ un entier
n\'egatif. On ne peut pas faire \`a partir de $C$ de biliaison \'el\'ementaire 
de hauteur $h$ sur $Q$.

\dem Quitte \`a tensoriser $\cE$ par un faisceau inversible, on peut
supposer qu'on  a  une suite exacte :
$$0\to {\cal O}_{\bf P}(-a)\to \cE \to {\cal J}_C \to 0$$ 
et qu'on a $H^0\cE(a-1)=0$, puisque $C$ est minimale pour
$\cE$. Remarquons qu'on a aussi $H^0{\cal
J}_C(a-1)=0$, donc
$s\geq a$ et $a+h-s$ est strictement n\'egatif.

On a aussi $\omega_C\simeq {\cal O}_C(a-4)$.

Supposons qu'il existe un  homomorphisme injectif $u_1 : {\cal J}_{C/Q}(-h)
\to {\cal O}_Q$. D'apr\`es 1.5 il correspond \`a un \'el\'ement non nul
$\theta_1$ de
$H^0\omega_C(4+h-s)=H^0 {\cal O}_C(a+h-s)$.

Soit $n_0=  \inf \{ n \in {{\bf Z}}\;\mid\; h^0\cO_C(n) \neq 0 \; \}$. D'apr\`es ce
qui pr\'ec\`ede, on a $n_0\leq a+h-s<0$ et $s-a+n_0\leq h<0$.

 Soit
$\theta_2$ un
\'el\'ement non nul  de
$H^0 {\cal O}_C(n_0)=H^0\omega_C(4-a+n_0)$  qui correspond d'apr\`es 1.5 \`a
un homomorphisme non nul $u_2 : {\cal J}_{C/Q} \to
{\cal O}_Q(s-a+n_0)$.

Pour $i\in \{1,2\}$ soient $C'_i$ et $C''_i$ les courbes associ\'ees \`a
$\theta_i$ comme on les a construites en 2.3. Alors $\theta_i=\epsilon_i
p_i$, o\`u
$p_i$ est la projection de  ${\cal O}_{C}$  sur un quotient ${\cal
O}_{C'_i}$ et
$\epsilon_1 $ (resp. $\epsilon_2 $) est une injection de ${\cal
O}_{C'_1}$ dans $ {\cal O}_{C}(a+h-s)$ (resp. de ${\cal
O}_{C'_2}$ dans $ {\cal O}_{C}(n_0)$) et $u_i$ se prolonge \`a
${\cal J}_{C'_i/Q}$. De plus, ${\cal J}_{C''_i}$  est
l'annulateur de ${\cal J}_{C'_i}/{\cal J}_C$.

Le produit $\theta_2 \theta_1$ est une section de ${\cal
O}_C(a+h-s+n_0)$ et il est nul par d\'efinition de $n_0$. On a donc
$\epsilon_2 p_2\epsilon_1 p_1=0$ (en fait il faudrait plut\^ot l'\'ecrire
$\epsilon_2 p_2\epsilon_1(-n_0) p_1(-n_0)=0$, mais on omettra les
d\'ecalages), et 
$ p_2\epsilon_1=0$ puisque $p_1$ est surjectif et $\epsilon_2$ injectif. On en
d\'eduit que
$C'_2$ est contenu dans le support du conoyau de $\epsilon_1$, qui est est
aussi le conoyau de $\theta_1$, donc dans la plus grande courbe contenue dans
ce support, c'est-\`a-dire
$C''_1$. On a donc des inclusions ${\cal J}_{C}\subseteq {\cal
J}_{C''_1}\subseteq {\cal J}_{C'_2}$, et puisque $u_2$ se prolonge \`a ${\cal
J}_{C'_2/Q}$, il  se prolonge
\'egalement \`a
${\cal J}_{C''_1/Q}$.

D'apr\`es 1.10, puisque $u_1$ est injectif, il en est de m\^eme de son
prolongement \`a 
${\cal J}_{C'_1/Q}$. On d\'efinit ainsi une biliaison \'el\'ementaire
descendante, de hauteur $h$ sur
$Q$, qui associe $\Gamma $ \`a $C$ et $\Gamma'_1$ \`a $C'_1$. On a donc un diagramme
commutatif de suites exactes :
$$\matrix {0&\to& {\cal J}_{C/Q}(-h)&\to &{\cal J}_{C'_1/Q}(-h)&\to &{\cal J}_{C'_1}/
{\cal J}_{C}(-h)&\to& 0 \cr
&&\mapdown {}&&\mapdown {}&&&& \cr
0&\to &{\cal J}_{\Gamma/Q}&\to &{\cal J}_{\Gamma'_1/Q}&\to &{\cal J}_{\Gamma'_1}/
{\cal J}_{\Gamma}&\to& 0 \cr}$$
dans lequel les deux fl\`eches verticales sont des isomorphismes induits par
$u_1$. On en d\'eduit que ${\cal J}_{\Gamma'_1}/
{\cal J}_{\Gamma}$ est isomorphe \`a ${\cal J}_{C'_1}/ {\cal
J}_{C}(-h)$ donc que l'annulateur de ${\cal J}_{\Gamma'_1}/
{\cal J}_{\Gamma}$ est \'egal \`a ${\cal J}_{C''_1}$. Cet annulateur contient
\'evidemment
${\cal J}_{\Gamma}$, donc ${\cal J}_{\Gamma}$ est contenu dans
${\cal J}_{C''_1}$.

D'autre part, $u_2$ se
prolonge en un homomorphisme non nul ${\cal J}_{C''_1/Q}\to {\cal O}_Q(s-a+n_0)$. Par
composition, on obtient un homomorphisme non nul (cf. 1.9) :
$${\cal
J}_{C/Q}(-h)\simeq {\cal J}_{\Gamma/Q}\to {\cal O}_Q(s-a+n_0)$$  donc par 1.5
un \'el\'ement non nul de $H^0 {\cal O}_C (n_0+h)$ ce qui donne une
contradiction.

\th {Th\'eor\`eme 2.5}. Soit  $C$ une courbe sous-canonique minimale pour
 un fibr\'e $\cE$ de rang 2 sur ${\bf P}^3$. Alors $C$ est minimale dans sa
classe de biliaison.

\dem
Soit $H_{\gamma,M}$ le sch\'ema de Hilbert des courbes \`a cohomologie et
module de Rao constants contenant $C$, qui est irr\'eductible (cf. [MDP1] V et
VI). Nous aurons besoin des r\'esultats des deux lemmes suivants :

\th {Lemme 2.6}. L'ensemble des courbes sous-canoniques de $H_{\gamma,M}$,
minimales pour un fibr\'e, est un ouvert.

\dem On a un isomorphisme $\omega_C\simeq {\cal O}_C(\alpha)$. Une courbe $C'$ de
$H_{\gamma,M}$ est sous-canonique si et seulement si il existe un entier
$\alpha'$ et un isomorphisme 
$\omega_{C'}\simeq {\cal O}_{C'}(\alpha')$. Mais puisque $C$ et $C'$ ont
m\^eme cohomologie, on a alors $\alpha=\alpha'$. L'isomorphisme $\omega_C\simeq {\cal
O}_C(\alpha)$ correspond \`a une section de $\omega_C(-\alpha)$ qui se prolonge sur un
voisinage $U$ de $C$ dans
$H_{\gamma,M}$ (cf. [MDP1] VII 2.3 et 2.5). Pour toute courbe $C'$ de $U$, il
existe un homomorphisme injectif  (quitte \`a restreindre $U$)
${\cal O}_{C'}\to \omega_{C'}(-\alpha)$ qui est un isomorphisme car les deux
faisceaux ont m\^eme polyn\^ome de Hilbert.

Soit $C'$ une courbe sous-canonique de $H_{\gamma,M}$. On voit facilement que
$C'$ est minimale pour le fibr\'e auquel elle correspond si et seulement si $H^0
{\cal J}_{C'}(\alpha+3)$ est nul et l'ensemble des courbes $C'$ de
$H_{\gamma,M}$ v\'erifiant $H^0
{\cal J}_{C'}(\alpha+3)=0$ est soit vide, soit \'egal \`a $H_{\gamma,M}$.

\th {Lemme 2.7}. L'ensemble des courbes de $H_{\gamma,M}$ pour lesquelles
on peut faire une biliaison $(s,h)$ est un ouvert.

\dem C'est une cons\'equence de [MDP1] VII 4.7.

\medbreak\noindent{\bf Fin de la d\'emonstration du Th\'eor\`eme 2.5.}\enspace
\nobreak Si
$C$ n'est pas minimale dans sa classe de biliaison, il existe un entier
$m\geq 1$, une suite de courbes 
$C_0,C_1, \ldots, C_m$ 
telle que $C_{i+1}$ s'obtienne \`a partir de $C_{i}$ par une biliaison
\'el\'ementaire de hauteur strictement positive, et $C$ \`a partir de $C_m$
par une d\'eformation
\`a cohomologie et module de Rao constants (cf [MDP1] IV 5). Soit $s$ le degr\'e
de la surface sur laquelle on fait la biliaison \'el\'ementaire qui fait
passer de $C_{m-1}$ \`a $C_m$ et $-h$ sa hauteur. L'ouvert des courbes de
$H_{\gamma,M}$ pour lesquelles on peut faire une biliaison $(s,h)$ est donc
non vide. Puisque $H_{\gamma,M}$ est irr\'eductible, cet ouvert rencontre
l'ouvert des courbes sous-canoniques minimales pour un fibr\'e,  et la
proposition 2.4 donne une contradiction.

\vskip 1cm
\titre {R\'ef\'erences bibliographiques}

[BBM] Ballico E., Bolondi G. Migliore J., The Lazarsfeld-Rao problem for
liaison classes of two-codimensional subschemes of $ {\bf P}^n$. Amer. J.
Math. 113,  117--128 (1991). 

[B] Buraggina A., Biliaison classes of reflexive
sheaves, Math. Nachr. 201, 53--76, 1999.

[Ho] Horrocks G., Vector bundles on the punctured spectrum of a local
ring, Proc. Lond. Math. Soc., 14, 689-713 (1964).

[HMDP]  Hartshorne R., Martin-Deschamps M. et  Perrin D., Un th\'eor\`eme
de Rao pour les familles de courbes gauches, Journal of Pure and Applied
Algebra Algebra 155, 53-76 (2001).

[MDP1]  Martin-Deschamps M. et 
Perrin D., Sur la classification des courbes gauches I, Ast\'erisque, Vol.
184-185, 1990.

[MDP2]  Martin-Deschamps M. et 
Perrin D.,  Quand un morphisme de fibr\'es d\'eg\'en\`ere-t-il le long d'une courbe
lisse ?
		Lecture Notes in Pure and Applied Mathematics Series/200. Marcel Dekker, Inc.
july 1998.

[Mi] Migliore J., Geometric Invariants of Liaison, J.  Algebra 99, 548-572
(1986).

[R] Rao A. P., Liaison among curves in $ {\bf P}^3$, Invent. Math., Vol.
50, 205.217 (1979).

\bye